\documentclass{article}
\usepackage[utf8]{inputenc}

\pdfoutput=1
\usepackage{mathtools}
\usepackage{amsfonts}
\usepackage{amssymb}
\usepackage{centernot}
\usepackage{graphicx}
\usepackage{amsthm}
\usepackage{enumerate}

\usepackage{tikz}
\usepackage{tikz-cd}
\usetikzlibrary{decorations.markings}
\usetikzlibrary{arrows}
\usetikzlibrary{calc}
\providecommand{\U}[1]{\protect\rule{.1in}{.1in}}
\newtheorem{theorem}{Theorem}

\newtheorem{definition}[theorem]{Definition}

\newtheorem{lemma}[theorem]{Lemma}

\theoremstyle{remark}

\newcommand{\Qr}{Q^{\textup{rev}}}
\newcommand{\br}{\beta^{\textup{rev}}}

\newcommand{\V}{\mathcal{V}}
\newcommand{\M}{\mathcal{M}}

\allowdisplaybreaks

\begin{document}

\title{Reorienting quandle orbits}
\author{Lorenzo Traldi\\Lafayette College\\Easton, PA 18042, USA\\traldil@lafayette.edu
}
\date{ }
\maketitle

\begin{abstract}
Motivated by knot theory, it is natural to define the orientation-reversal of a quandle orbit by inverting all the translations given by elements of that orbit. In this short note we observe that this natural notion is unsuited to medial quandles.

\end{abstract}

\section{Introduction}

A \emph{quandle} is defined by a binary operation $\triangleright$ on a set $Q$, which satisfies these three axioms.

\begin{itemize}
\item $ x \triangleright x = x \thickspace \forall x \in Q.$
\item $ \textup{For each }y \in Q, \textup{ the formula }\beta_y(x)=x \triangleright y \textup{ defines a bijection }  \beta_y:Q \to Q.$
\item $ (x \triangleright y) \triangleright z=(x \triangleright z) \triangleright (y \triangleright z)  \thickspace \forall x,y,z \in Q.$ 
\end{itemize}

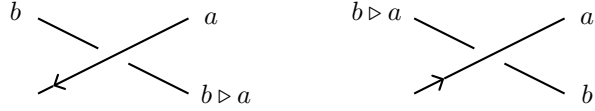
\begin{figure} [t]
\centering 
\begin{tikzpicture} [>=angle 90]
\draw [thick] [->] (1,.5) -- (-0.8,-0.4);
\draw [thick] (-0.8,-0.4) -- (-1,-.5);
\draw [thick] (-1,.5) -- (-.2,0.1);
\draw [thick] (0.2,-0.1) -- (1,-.5);
\node at (1.3,0.5) {$a$};
\node at (-1.3,0.6) {$b$};
\node at (1.5,-0.5) {$b \triangleright a$};

\draw [thick]  (6,.5) -- (4.4,-0.3);
\draw [thick] [<-] (4.4,-0.3) -- (4,-.5);
\draw [thick] (4,.5) -- (4.8,0.1);
\draw [thick] (5.2,-0.1) -- (6,-.5);
\node at (6.3,0.5) {$a$};
\node at (3.5,0.6) {$b \triangleright a$};
\node at (6.3,-0.5) {$b$};
\end{tikzpicture}
\caption{The underpassing arc $b$ is on the right of $a$, and $b \triangleright a$ is on the left.}
\label{fig1}
\end{figure}

Special cases of this definition were considered almost 100 years ago, but the modern theory dates back to the introduction of quandles in classical knot theory by Joyce \cite{J} and Matveev \cite{M} in the 1980s. The central object of classical knot theory is a \emph{link}, i.e., a collection of finitely many  simple closed curves in $\mathbb R ^3$; the curves are oriented, pairwise disjoint, and piecewise smooth. Each individual curve is a \emph{component} of the link. Links are represented by diagrams, and the fundamental quandle of a link is generated by the arcs of a diagram. The operation $\triangleright$ of this quandle may be interpreted as ``passing under, from right to left.'' See Figure \ref{fig1}.

Recall that a quandle variety $\mathcal V$ is the collection of quandles that satisfy some particular set of formulas involving $\triangleright$ and $\triangleright^{-1}$. Two well-known examples are the variety $\M$ of \emph{medial} quandles, defined by $(w \triangleright x) \triangleright (y \triangleright z)= (w \triangleright y) \triangleright (x \triangleright z) \thickspace \forall w,x,y,z \in Q$, and the variety $\V_n$ of $n$-quandles, defined by $\beta_y^n(x) =x \thickspace \forall x,y \in Q$. (Here $n$ is an integer.)  For thorough discussions of the general theory of medial quandles we refer to Jedli\v{c}ka \emph{et al.} \cite{JPSZ1, JPSZ2}. Note that some authors, like Joyce \cite{J}, call medial quandles ``abelian.''

If $Q$ is an arbitrary quandle and $\mathcal V$ is a quandle variety, then $Q$ gives rise to an element $Q_{\mathcal V} \in \mathcal V$ in a natural way: The formulas that define $\mathcal V$ yield an equivalence relation on $Q$, and $Q_{\mathcal V}$ is the quotient of $Q$ modulo this equivalence relation. In particular, Joyce \cite[Sec.\ 10]{J} refers to $Q_{\M}$ as the ``largest abelian quotient of $Q$.''

A basic operation of knot theory is that we may reverse the orientation of one of the components of a link. As is clear from Figure \ref{fig1}, the effect of this operation on the link quandle is to interchange the functions $\beta_a$ and $\beta_a^{-1}$ for each arc $a$ that belongs to the reversed component. 

If $Q$ is a quandle and $x \in Q$ then $Q_x$ denotes the orbit of $x$ in $Q$, i.e., the smallest subset of $Q$ that contains $x$ and is closed under $\beta_y$ and $\beta_y^{-1}$ for every $y \in Q$. It is not hard to see that in a link quandle, the diagram arcs corresponding to a single component all belong to the same orbit. If that component's orientation is reversed, then every one of these arcs has its orientation reversed. These properties suggest the following.

\begin{definition}
\label{orrev}
Let $Q$ be a quandle with an orbit $Q_x$. Then the quandle obtained from $Q$ by \emph{reversing the orientation} of $Q_x$ is the quandle $\Qr(x)$ with the same underlying set as $Q$, in which $\br_y=\beta_y\thickspace \allowbreak \forall y \notin Q_x$ and $\br_y=\beta_y ^{-1} \thickspace \allowbreak \forall y \in Q_x$.
\end{definition}

Verifying that Definition \ref{orrev} really does yield a quandle $\Qr(x)$ is not difficult; we leave the details to the reader. We should say that although we have not seen Definition \ref{orrev} in the literature, we presume that such an easy and intuitive idea has occurred to other researchers.

Our purpose in writing this note is to bring up the problem of defining a notion of orientation-reversal that is appropriate for quandles in a variety $\V$. For $n$-quandles there is no difficulty: $\beta_y^n$ is the identity map if and only if $\beta_y^{-n}$ is the identity map, so $Q$ is an $n$-quandle if and only if $\Qr(x)$ is an $n$-quandle. But for other varieties, there is no reason to expect the defining formulas to be satisfied when some $\beta$ maps are replaced by their inverses. For such a variety, the natural way to define orientation-reversal is to consider the quotient $\Qr(x)_{\V}$, which may be a smaller quandle than $\Qr(x)$.

Observe that if we reverse the orientation of the same orbit twice, the result is the original quandle: $(\Qr(x))^{\textup{rev}}(x)=Q$. It follows that $\Qr(x)$ determines $Q$. This observation motivates our next definition.

\begin{definition}
\label{suit}
A quandle variety $\V$ is \emph{suited for orientation-reversals} if it has the property that for every $Q \in \V$, $\Qr(x)_{\mathcal{V}}$ determines $Q$.
\end{definition}

In the remainder of the note we show that $\M$ is not suited for orientation-reversals, by providing an example of an infinite medial quandle $Q$ with the property that $\Qr(x)_{\M}$ is a trivial quandle with only two elements.

\section{An example}

Let $\Lambda = \mathbb Z[t,t^{-1}]$ be the ring of Laurent polynomials in the variable $t$, with integer coefficients. Let $J$ be the principal ideal of $\Lambda$ generated by $t^2+t-1$, and let $N=\Lambda/J$. Notice that if $x\in N$ then $(1-t)x=(1-t)x+(t^2+t-1)x=t^2 x$ and hence $x=(t^2+t)x$. We leave it to the reader to use these properties to prove the following.

\begin{lemma}
\label{element}
Every element of $N$ is $x=n_1t+n_2t^2+J$ for some unique $n_1,n_2 \in \mathbb Z$.
\end{lemma}

In the terminology of Jedli\v{c}ka \emph{et al.} \cite{JPZ}, our example $Q$ is a free $(t^2+t-1)$-quandle on a two-element set. It is described explicitly as follows. Let $Q_1$ and $Q_2$ be two disjoint copies of $N$. If $x \in N$, we use $x_i$ to denote the copy of $x$ in $Q_i$. Also, let $m_1=0+J$ and $m_2=1+J$. We build a medial quandle $Q$ on $Q_1 \cup Q_2$ by defining
\[
\beta_{y_j}(x_i)=x_i \triangleright y_j = (m_j-m_i+tx + (1-t)y)_i \in Q_i.
\]
Elementary arguments verify that $Q$ is a medial quandle, and
\[
\beta_{y_j}^{-1}(x_i)=t^{-1} \cdot  (m_i-m_j+x - (1-t)y)_i \in Q_i.
\]

\begin{lemma}
\label{orbits}
The orbits of $Q$ are $Q_1$ and $Q_2$.
\end{lemma}
\begin{proof}
If $x \in N$ then 
\[
0_1 \triangleright (t^{-2}(x-1))_2 = (1+t \cdot 0+(1-t)(t^{-2}(x-1)))_1 
\]
\[
= (1+t \cdot 0+t^2(t^{-2}(x-1)))_1 =(1+x-1)_1=x_1 \text{,}
\]
so $x_1$ is in the same orbit of $Q$ as $0_1$. Similarly, $Q_{x_2}=Q_{0_2}$ because
\[
0_2 \triangleright (t^{-2}(x+1))_1 = (-1+t \cdot 0+(1-t)(t^{-2}(x+1)))_2 
\]
\[
= (-1+t \cdot 0+t^2(t^{-2}(x+1)))_2 =(-1+x+1)_2=x_2 \text{.} 
\]

\end{proof}

Notice that for convenience, in equations we sometimes suppress the coset notation for an element of $N$. For instance, in the last line of the proof of Lemma \ref{orbits} we write $(-1+x+1)_2$ rather than $(-1+J+x+1+J)_2$.

\begin{theorem}
$\Qr(0_2)_\M$ has only two elements.
\end{theorem}
\begin{proof}
Suppose $w,x,y,z \in N$. Using the formulas for $\beta_{y_j}$ and $\beta^{-1}_{y_j}$ displayed above, in $\Qr(0_2)$ we have
\begin{align*}
(w_1 \triangleright x_1) \triangleright { } & (y_1 \triangleright z_2) = (tw+(1-t)x)_1 \triangleright (t^{-1} \cdot (m_1-m_2+y-(1-t)z))_1
\\
&=(t(tw+(1-t)x)+(1-t)(t^{-1} \cdot (-1+y-(1-t)z)))_1
\\
&=(t^2w+t(1-t)x+t^2(t^{-1} \cdot (-1+y-t^2 z)))_1
\\
&=(t^2w+t^3x+t(-1+y-t^2 z))_1=(t^2w+t^3x-t+ty-t^3 z)_1.
\end{align*}
In $\Qr(0_2)_{\M}$ we have $(w_1 \triangleright x_1) \triangleright (y_1 \triangleright z_2)=(w_1 \triangleright y_1) \triangleright (x_1 \triangleright z_2)$, and therefore
\begin{equation}
\label{require}
(t^2w+t^3x-t+ty-t^3 z)_1=(t^2w+t^3y-t+tx-t^3 z)_1.
\end{equation}

Let $a \in N$ be arbitrary. Taking $w=x=1+J,y=0+J$ and $z=-t^{-3}a$ in (\ref{require}) yields
$(t^2+t^3-t+a)_1=(t^2-t+t+a)_1$. As $t^2+t^3-t=t(t^2+t-1) \in J$, it follows that $a_1 = (t^2+a)_1$ in $\Qr(0_2)_{\M}$. Taking $w=x=0+J,y=t^{-2}+J$ and $z=-t^{-3}(a-1)$ in (\ref{require}) yields $(-t+t^{-1}+a-1)_1=(t-t+a-1)_1$. As $-t+t^{-1}-1=-t^{-1}(t^2+t-1) \in J$, it follows that $a_1=(a-1)_1$ in $\Qr(0_2)_{\M}$. Using the facts that $a_1 = (t^2+a)_1$ and $a_1=(a-1)_1$  in $\Qr(0_2)_{\M}$ for all $a \in N$, we conclude that also $a_1=(a+t^2+t-1)_1 = (a+t-1)_1=(a+t)_1$ in $\Qr(0_2)_{\M}$ for all $a \in N$. According to Lemma \ref{element}, the equalities $a_1=(a+t^2)_1$ and $a_1=(a+t)_1$ imply that $a_1=0_1$ in $\Qr(0_2)_{\M}$ for all $a \in N$. That is, the orbit of $0_1$ in $\Qr(0_2)_{\M}$ has only one element.

The proof of Lemma \ref{orbits} shows that every element of the orbit $Q_{0_2}$ is $0_2 \triangleright y_1$ for some element $y_1$ of the orbit $Q_{0_1}$. The functions $\beta_{y_1}$ are the same in $Q$ and $\Qr(0_2)$, so it is also true that every element of $\Qr(0_2)_{0_2}$ is $0_2 \triangleright y_1$ for some $y_1 \in \Qr(0_2)_{0_1}$. Of course this property is inherited by the quotient quandle $\Qr(0_2)_{\M}$. But in the previous paragraph we showed that $(\Qr(0_2)_{\M}) _{0_1}$ has only one element. Therefore $(\Qr(0_2)_{\M}) _{0_2}$ also has only one element.
\end{proof}

\end{document}